\documentclass[journal]{IEEEtran}
\IEEEoverridecommandlockouts           
\usepackage{graphics} 
\usepackage{epsfig} 
\usepackage{mathptmx} 
\usepackage{times} 
\usepackage{amsmath} 
\usepackage{amssymb}  
\newcommand{\R}{{\mathbb R}}

\usepackage{setspace,verbatim,longtable,dcolumn,ifpdf}
\usepackage{cite}
\usepackage{color}

\usepackage{amsfonts}

\newtheorem{cor}{Corollary}

\newtheorem{theorem}{Theorem}
\newtheorem{definition}[theorem]{Definition}
\newtheorem{lem}[theorem]{Lemma}

\providecommand{\norm}[1]{\left \| #1\right \|}

\renewcommand{\S}{{\mathbb S}}

\newcommand{\bmat}[1]{\begin{bmatrix}#1\end{bmatrix}}

\begin{document}

\title{A Converse Sum of Squares Lyapunov Result with a Degree Bound}
\author{Matthew~M.~Peet and Antonis~Papachristodoulou
\thanks{M. M. Peet is with the Department of Mechanical, Materials, and Aerospace Engineering, Illinois Institute of Technology, 10 West 32nd Street, E1-252B, Chicago, IL 60616, U.S.A. {\tt\small mpeet@iit.edu}}
\thanks{A.~Papachristodoulou is with the Department of Engineering Science, University of Oxford, Parks Road, Oxford, OX1 3PJ, U.K. {\tt\small antonis@eng.ox.ac.uk}.}%
\thanks{This material is based upon work supported by the National Science Foundation under Grant No. CMMI 110036 and from
EPSRC grants EP/H03062X/1, EP/I031944/1 and EP/J012041/1}}
\maketitle

\begin{abstract}
Sum of Squares programming has been used extensively over the past decade for the stability analysis of nonlinear systems but several questions remain unanswered. In this paper, we show that exponential stability of a polynomial vector field on a bounded set implies the existence of a Lyapunov function which is a sum-of-squares of polynomials. In particular, the main result states that if a system is exponentially stable on a bounded nonempty set, then there exists an SOS Lyapunov function which is exponentially decreasing on that bounded set. The proof is constructive and uses the Picard iteration. A bound on the degree of this converse Lyapunov function is also given. This result implies that semidefinite programming can be used to answer the question of stability of a polynomial vector field with a bound on complexity.
\end{abstract}

%

\section{INTRODUCTION}


Computational and numerical algorithms are extensively used in control theory. A particular example is semidefinite programming conditions for addressing linear control problems, which are formulated as Linear Matrix Inequalities (LMIs). Using such tools, several questions on the analysis and synthesis of linear systems can be formulated and addressed effectively. In fact, ever since the 1990s~\cite{BoyEFB94}, LMIs have had a significant impact in the control field, to the point that once the solution of a control problem has been formulated as the solution to an LMI, it is considered solved.

When it comes to nonlinear and infinite-dimensional systems, the equivalent problems can be formulated as polynomial non-negativity constraints under a Lyapunov framework, but these are not, at first glance, as easy to solve. Polynomial non-negativity is in fact NP hard. It is for this reason that several researchers have looked at alternative tests for non-negativity, that are polynomial-time complex to test, and which imply non-negativity. One such relaxation is the existence of a sum of squares decomposition: the ability to optimize over the set of positive polynomials using the sum-of-squares relaxation has undoubtedly opened up new ways for addressing nonlinear control problems, in much the same way Linear Matrix Inequalities are used to address analysis questions for linear finite-dimensional systems. However, there remain several open questions about how these methods can be used to search for Lyapunov functions for nonlinear systems. For references on early work on optimization of polynomials, see~\cite{lasserre_2001,nesterov_2000}, and~\cite{parrilo_thesis}. For more recent work see~\cite{henrion_2005} and~\cite{chesi_2007}. For a recent review paper, see~\cite{Che10}. Today, there exist a number of software packages for optimization over positive polynomials, e.g. SOSTOOLS~\cite{prajna_2004} and GloptiPoly~\cite{henrion_2001}.

At the same time, there are still a number of unanswered questions regarding the use of sum of squares as a relaxation to nonnegativity and its use for the analysis of nonlinear systems. Unanswered questions include, for example, a series of questions on controller synthesis and the role of duality to convexify this problem, as well as estimating regions of attraction of equilibria. On the computation and optimization side, it is unclear whether multi-core computing could be used for computation, as well as how to take advantage of sparsity in semidefinite programming.

In this paper, we do not consider the problem of computing sum-of-squares Lyapunov functions. Such work can be found in, e.g.~\cite{parrilo_thesis,Papachristodoulou_2002,wang_thesis,tan_2006}. Instead, we concentrate on the properties of the converse Lyapunov functions for systems of the form
\begin{equation*}
    \dot{x}(t)=f(x(t)),
\end{equation*}
where $f:\R^n \rightarrow \R^n$ is polynomial. In particular, we address the question of whether an exponentially stable nonlinear system will have a sum-of-squares Lyapunov function which establishes this property. This result adds to our previous work~\cite{peet_2009a}, where we were able to show that exponential stability on a bounded set implies the existence of a exponentially decreasing polynomial Lyapunov function on that set.

Work that is relevant to the one presented here includes research on continuity properties, see e.g.~\cite{barbasin_1951},~\cite{malkin_1954} and~\cite{kurzweil_1963} and the overview in~\cite{massera_1956}. Infinitely-differentiable functions were explored in the work~\cite{wilson_1969,lin_1996}. Other innovative results are found in~\cite{lakshmikantam_1992} and~\cite{teel_1999}. The books~\cite{hahn_1967} and~\cite{krasovskii_1963} treat further converse theorems of Lyapunov. Continuity of Lyapunov functions is inherited from continuity of the solution map with respect to initial condition. An excellent treatment of this problem can be found in the text of Arnol'd~\cite{arnold_book}.

Unlike the work in~\cite{peet_2009a}, this paper is closely tied to systems theory as opposed to approximation theory. Our method is to take a well-known form of converse Lyapunov function based on the solution map and use the Picard iteration to approximate the solution map. The advantage of this approach is that if the vector field is polynomial, the Picard iteration will also be polynomial. Furthermore, the Picard iteration inductively retains almost all the properties of the solution map. The result is a new form of iterative converse Lyapunov function, $V_k$. This function is discussed in Section~\ref{se:mainresult}.

The first practical contribution of this paper is to give a bound on the number of decision variables involved in the question of exponential stability of polynomial vector fields on bounded sets. This is because SOS functions of bounded degree can be parameterized by the set of positive matrices of fixed size. Furthermore, we note that the question of existence of a Lyapunov function with negative derivative is convex. Therefore, if the question of polynomial positivity on a bounded set is decidable, we can conclude that the problem of exponential stability of polynomial vector fields on that set is decidable. The further complexity benefit of using SOS Lyapunov functions is discussed in Section~\ref{se:benefits}.

%
%

The main result of the paper is stated and proven in Section~\ref{se:mainresult}. Preceding the main result is a series of lemmas that are used in the proof of the main theorem. In Subsection~\ref{se:picarditeration} we show that the Picard iteration is contractive on a certain metric space; and in Subsection~\ref{se:picardextensionconvergence} we propose a new way of extending the Picard iteration. In Section~\ref{se:derivativeinequalitylemma} we show that the Picard iteration approximately retains the differentiability properties of the solution map, before we prove the main result. The implications of the main result are then explored in Section~\ref{se:benefits} and Section~\ref{se:quadratic}. A detailed example is given in Section~\ref{se:example}. The paper is concluded in Section~\ref{se:conclusion}.

\section{Main Result}

Before we begin the technical part of the paper, we give a simplified version of the main result.
\begin{theorem}
Suppose that $f$ is polynomial of degree $q$ and that solutions of $\dot x = f(x)$ satisfy
\begin{equation*}
\norm{x(t)}\le K \norm{x(0)}e^{-\lambda t}
\end{equation*}
for some $\lambda>0$, $K\ge1$ and for any $x(0)\in M$, where $M$ is a bounded nonempty region of radius $r$. Then there exist $\alpha,\beta,\gamma>0$ and a sum-of-squares polynomial $V(x)$ such that for any $x \in M$,
\begin{align*}
\alpha \norm{x}^{2} \le V(x)\le \beta \norm{x}^{2}  \\
\nabla V(x)^T f(x)\le  - \gamma \norm{x}^{2}.
\end{align*}
Further, the degree of $V$ will be less than $2q^{(N k-1)}$,
where $k(L,\lambda,K)$ is any integer such that
$c(k):=\sum_{i=0}^{N-1} \left(e^{TL} + K (TL)^k\right)^i K^2 (TL)^k <K$, and
\begin{align*}
&c(k)^2 + \frac{\log 2K^2}{2 \lambda} K \frac{(TL)^k}{T} (1+c(k)) (K+c(k)) < \frac{1}{2},  \\
&c(k)^{2} <  \frac{ \lambda}{KL\log 2K^2}(1-(2K^2)^{-\frac{L}{\lambda}})
\end{align*}
  and $N(L,\lambda,K)$ is any integer such that $NT>\frac{\log 2K^2}{2 \lambda}$ and $T<\frac{1}{2L}$ for some $T$ and where $L$ is a Lipschitz bound on $f$ on $B_{4Kr}$.
\end{theorem}

\section{Sum-of-Squares}\label{sec:SOS}

Sum of squares (SOS) methods have been introduced over the past decade to allow for the algorithmic solution of problems that frequently arise in systems and control theory, many of which can be formulated as polynomial non-negativity constraints that are however difficult to solve. In these methods, non-negativity is relaxed to the existence of a SOS decomposition, which can be tested using Semidefinite programming.

Consider, for example, the problem of ensuring that a polynomial $p(x)\in \R[x]$ satisfies $p(x) \geq 0$. This problem arises naturally when trying to construct Lyapunov functions for the stability analysis of dynamical systems, which is the topic of this paper. Since ensuring non-negativity is hard~\cite{MurK87} many researchers have investigated alternative ways to do this. In~\cite{Sho87}, the existence of a Sum of Squares decomposition was used for that purpose, which involves the presentation of other polynomials $p_i(x)$ such that
\begin{equation}
  p(x) = \sum_{i = 1}^k p_i(x)^2 \label{eq:SOS}
\end{equation}
Algorithms for ensuring this have appeared in the 1990's~\cite{PowW98} but it was not until the turn of the century that this was recognized as being solvable using Semidefinite Programming~\cite{Par00}. In particular, (\ref{eq:SOS}) can be shown equivalent to the existence of a $Q  \succeq 0$ and a vector of monomials $Z(x)$ of degree less than or equal half the degree of $p(x)$, such that
\begin{equation*}
  p(x) = Z(x)^T Q Z(x)
\end{equation*}
In the above representation, the matrix $Q$ is not unique, in fact it can be represented as
\begin{equation}
  Q = Q_0 + \sum_i \lambda_i Q_i \label{eq:SOSLMI}
\end{equation}
where $Q_i$ satisfy $Z(x)^T Q_i Z(x) = 0$. The search for $\lambda_i$ such that $Q$ in (\ref{eq:SOSLMI}) is such that $Q  \succeq 0$ is a Linear Matrix Inequality, which can be solved using Semidefinite Programming. Moreover, if $p(x)$ has unknown coefficients that enter affinely in the representation (\ref{eq:SOS}), Semidefinite Programming can be used to find values for them so that the resulting polynomial is SOS.

This latter observation can allow us to \emph{search} for polynomials that satisfy SOS conditions: the most important example is in the construction of Lyapunov functions, which is the topic of this paper. For more details, please see~\cite{Par00,Papachristodoulou_2002}. The question that we address in this paper is whether Sum of Squares Lyapunov functions always exist for locally exponentially stable systems.

\section{Notation and Background}

The core concept we use in this paper is the Picard iteration. We use this to construct an approximation to the solution map and then use the approximate solution map to construct the Lyapunov function. Construction of the Lyapunov function will be discussed in more depth later on.

Denote the Euclidean ball centered at $0$ of radius $r$ by $B_r$. Consider an ordinary differential equation of the form
\begin{equation}
\dot{x}(t)=f(x(t)), \qquad x(0)=x_0,\qquad f(0)=0, \label{eq:sys1}
\end{equation}
where $x \in \R^n$ and $f$ satisfies appropriate smoothness properties for local existence and uniqueness of solutions. The solution map is a function $\phi$ which satisfies
\begin{equation*}
\frac{\partial}{\partial t}\phi(t,x)=f(\phi(t,x)) \qquad \text{ and }\qquad \phi(0,x)=x.
\end{equation*}

\subsection{Lyapunov Stability}

The use of Lyapunov functions to prove stability of ordinary differential equations is well-established. The following theorem illustrates the use of Lyapunov functions.
\begin{definition}
We say that the system defined by the equations in~\eqref{eq:sys1} are exponentially stable on $X$ if there exist $\gamma, K >0$ such that for any $x_0 \in X$,
\[
\norm{x(t)}\le K \norm{x_0}e^{-\gamma t}
\]
for all $t \ge 0$.
\end{definition}

\begin{theorem}[Lyapunov]
Suppose there exist constants $\alpha,\beta,\gamma>0$ and a continuously differentiable function $V$ such that the following conditions are satisfied for all $x \in U \subset \R^n$.
\begin{align*}
\alpha \norm{x}^2 \le V(x)  \le \beta \norm{x}^2\\
\nabla V(x)^T f(x) \le -\gamma \norm{x}^2
\end{align*}
Then we have exponential stability of System~\eqref{eq:sys1} on $\left\{ x \;:\; \{y: V(y) \le V(x)\}\subset U\right\}$.
\end{theorem}


\subsection{Fixed-Point Theorems}

\begin{definition}
Let $X$ be a metric space. A mapping $F:X\rightarrow X$ is \textit{contractive} with coefficient $d \in [0,1)$ if
\begin{equation*}
\norm{Fx-Fy}\le d \norm{x-y} \qquad x,y \in X.
\end{equation*}
\end{definition}

The following is a \textit{Fixed-Point} Theorem.

\begin{theorem}[Contraction Mapping Principle~\cite{marsden_analysis}] Let $X$ be a complete metric space and let $F:X\rightarrow X$ be a contraction with coefficient $d$. Then there exists a unique $y \in X$ such that
\begin{equation*}F y=y.
\end{equation*}
Furthermore, for any $x_0 \in X$,
\begin{equation*}
\norm{F^k x_0 - y}\le d^k \norm{x_0-y}.
\end{equation*}
\end{theorem}

To apply these results to the existence of the solution map, we use the Picard iteration.
%

\section{Picard Iteration}\label{se:picarditeration}

We begin by reviewing the Picard iteration. This is the basic mathematical tool we will use
to define our approximation to the solution map and can be found in many texts, e.g.\cite{coddingtion_1955}.

\begin{definition}For given $T$ and $r$, define the complete metric space
\begin{align}
X_{T,r}:=\left\{q(t) : \begin{array}{l} \sup_{\substack{t\in[0,T]}} \norm{q(t)} \le r, \text{ }
q \text{ is continuous}. \end{array} \right\}
\end{align}
with norm
\[
\displaystyle \norm{q}_X=\sup_{\substack{t\in[0,T]}}\norm{q(t)}.
\]
\end{definition}

For a fixed $x \in B_r$ and $q \in X_{T,r}$, the \textit{Picard Iteration}~\cite{lindelhof_1894}, is defined as
\begin{equation*}
(P q)(t) \triangleq x + \int_{0}^t f(q(s))ds.
\end{equation*}
In this paper, we also define the Picard iteration iteration on functions $z(t,x)$ as
\begin{equation*}
(P z)(x,t) \triangleq x + \int_{0}^t f(z(x,s))ds.
\end{equation*}

We begin by showing that for any radius $r$, there exists a $T$ such that the Picard iteration is contractive on $X_{T,2r}$ for any $x \in B_r$.

\begin{lem}\label{lem:contraction_Pk} Given $r>0$, let $T<\min\{\frac{r}{Q},\frac{1}{L}\}$ where $f$ has Lipschitz factor $L$ on $B_{2r}$ and $Q=\sup_{x \in B_{2r}}f(x)$. Then $P:X_{T,2r} \rightarrow X_{T,2r}$ and there exists some $\phi \in X_{T,2r}$ such that for $t \in [0,T]$ and $x \in B_r$,
\begin{equation*}
\frac{d}{dt}\phi(t)=f(\phi(t)),\qquad \phi(0)=x
\end{equation*}
and for any $z \in X_{T,2r}$,
\begin{equation*}
\norm{\phi - P^k z }\le (TL)^k \norm{\phi-z}.
\end{equation*}
\end{lem}
\begin{proof}
We first show that for $x \in B_r$, $P:X_{T,2r} \rightarrow X_{T,2r}$. If $q \in X_{T,2r}$, then $\sup_{\substack{t\in[0,T]}} \norm{q(t)}\le 2r$ and so
\begin{align*}
\norm{Pq}_X &= \sup_{\substack{t\in[0,T]}}  \norm{ x + \int_0^{t} f( q(s) )}ds\\
&\le \norm{x}+\int_0^{T} \norm{f(q(s))}ds\\
&\le r + \int_0^{T} \sup_{y \in B_{2r}} \norm{f(y)}ds\\
& \le r + T Q < 2r
\end{align*}
Thus we conclude that $Pq \in X_{T,2r}$. Furthermore, for $q_1,q_2 \in X_{T,2r}$,
\begin{align*}
\norm{Pq_1 - Pq_2}_X&=\sup_{t \in [0,T]}\left\| \int_{0}^t\left(f(q_1(s))-f(q_2(s))\right)ds \right\|\\
&\le \int_{0}^{T} \sup_{\substack{t\in[0,T]}} \left\| f(q_1(s))-f(q_2(s)) \right\|ds\\
&\le T L \sup_{\substack{t\in[0,T]}} \left\| q_1(s)-q_2(s) \right\| = T L \left\| q_1-q_2\right\|_X
\end{align*}

Therefore, by the contraction mapping theorem, the Picard iteration converges on $[0,T]$ with convergence rate $(TL)^k$.
\end{proof}

%
%

\subsection{Picard Extension Convergence Lemma} \label{se:picardextensionconvergence}

\begin{figure}
\centering
  \includegraphics[scale=.8]{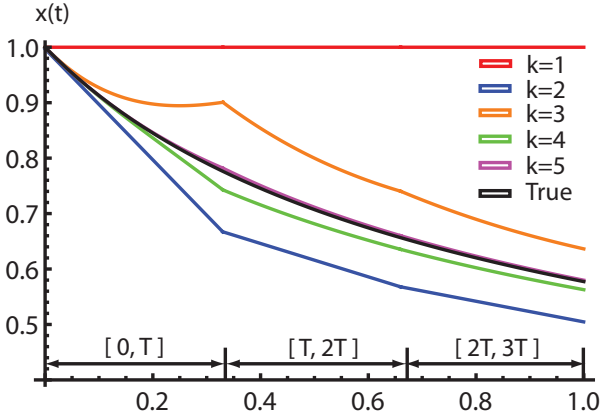}
 \caption{The Solution map $\phi$ and the functions $G^k_{i}$ for $k=1,2,3,4,5$ and $i=1,2,3$ for the system $\dot x(t) = - x(t)^3$. The interval of convergence of the Picard Iteration is $T=\frac{1}{3}$.}\label{fig:illustration}
\end{figure}

In this section we propose an extension to the Picard iteration approximation. We divide the interval into subintervals on which the Picard iteration is guaranteed to converge. On each interval, we apply the Picard iteration using the final value of the solution estimate from the previous interval as the initial condition, $x$. For a polynomial vector field, the result is a piece-wise polynomial approximation which is guaranteed to converge on an arbitrary interval -- see Figure~\ref{fig:illustration} for an illustration.

\begin{definition}
Suppose that the solution map $\phi$ exists on $ t \in [0,\infty)$ and $\norm{\phi(t,x)}\le K \norm{x}$ for any $x \in B_{r}$. Suppose that $f$ has Lipschitz factor $L$ on $B_{4Kr}$ and is bounded on $B_{4Kr}$ with bound $Q$. Given $T<\min\{\frac{2Kr}{Q},\frac{1}{L}\}$, let $z=0$ and define
\begin{equation*}
G^k_0(t,x):=(P^k z)(t,x)
\end{equation*}
and for $i>0$, define the functions $G_i$ recursively as
\begin{equation*}
G^k_{i+1}(t,x):=(P^k z)(t,G^k_i(T,x)).
\end{equation*}
The $G^k_i$ are Picard iterations $P^k z(t,x)$ defined on each sub-interval where we substitute the initial condition $x \mapsto G^k_{i-1}(t,x)$. Define the concatenation of these $G^k_i$ as
\begin{equation*}
G^k(t,x):=G^k_i(t-iT,x) \quad \forall \quad t\in [iT,iT+T] \quad \text{and } i=1,\cdots,\infty. 
\end{equation*}
\end{definition}

If $f$ is polynomial, then the $G^k_i$ are polynomial for any $i,k$ and $G^k$ is continuously differentiable in $x$ for any $k$.
The following lemma provides several properties for the functions $G^k$.

\begin{lem}\label{lem:converge_Gk}
Given $\delta>0$, suppose that the solution map $\phi(t,x)$ exists on $t \in [0,\delta]$ and on $x\in B_{r}$. Further suppose that $\norm{\phi(t,a,x)}\le K \norm{x}$ for any $x \in B_r$. Suppose that $f$ is Lipschitz on $B_{4Kr}$ with factor $L$ and bounded with bound $Q$. Choose $T<\min\{\frac{2Kr}{Q},\frac{1}{L}\}$ and integer $N>\delta/T$. Then let $G^k$ and $G_i^k$ be defined as above.

Define the function
\begin{equation*}
c(k)=\sum_{i=0}^{N-1} \left(e^{TL} + K (TL)^k\right)^i K^2 (TL)^k.
\end{equation*}
Given any $k$ sufficiently large so that $c(k)<K$, then
%
for any $x \in B_r$,
\begin{equation}
\sup_{s\in [0,\delta]}\norm{G^k(s,x)-\phi(s,x)}\le  c(k)\norm{x}.\label{eqn:CDN2}
\end{equation}
\end{lem}
\begin{proof}
Suppose $x \in B_r$. By assumption, the conditions of Lemma~\ref{lem:contraction_Pk} are satisfied using $r'=2Kr$. Let $z(t,x)=0$. Define the convergence rate $d=TL<1$. By Lemma~\ref{lem:contraction_Pk},
\begin{align*}
\sup_{s\in [0,T]}\norm{G_0^k(s,x)-\phi(s,x)} &=\sup_{\substack{s\in [0,T]}}\norm{(P^k z) (s,x)-\phi(s,x)}\\
 &\le d^k \sup_{s\in [0,T]}\norm{\phi(s,x)}\le K d^k\norm{x}.
\end{align*}
Thus Equation~\eqref{eqn:CDN2} is satisfied on the interval $[0,T]$. We proceed by induction. Define
\begin{equation*}
c_i(k)=\sum_{j=1}^i (e^{d} + K d^k)^j K^2 d^k.
\end{equation*}
and suppose that $\norm{G^k-\phi}_\infty\le c_{i-1}(k)\norm{x}$ on interval $[iT-T,iT]$. Then
\begin{align*}
&\sup_{s \in [iT,iT+T]}\norm{G^k(s,x)-\phi(s,x)}\\
& =\sup_{s \in [iT,iT+T]}\norm{G^k_i(s-iT,x)-\phi(s,x)}\\%
&=\sup_{s \in [iT,iT+T]}\norm{P^k z (s-iT,G_{i-1}^k(T,x))-\phi(s-iT,\phi(iT,x))}\\%
& \le \hspace{-1mm} \sup_{s \in [iT,iT+T]}\norm{P^k z (s-iT,G_{i-1}^k(T,x))-\phi(s-iT,G_{i-1}^k(T,x))}\\
& +\sup_{s \in [iT,iT+T]} \norm{\phi(s-iT,G_{i-1}^k(T,x))-\phi(s-iT,\phi(iT,x))}
\end{align*}
We treat these final two terms separately. First note that
\begin{align*}
\norm{G_{i-1}^k(T,x)} &\le \norm{\phi(iT,x)}+\norm{\phi(iT,x)-G_{i-1}^k(T,x)}\\
& \le K\norm{x} + c_{i-1}(k) \norm{x}\\
&\le (K+c_{i-1}(k))\norm{x}.
\end{align*}
Since $c_{i-1}(k)\le c(k)<K$ and $x \in B_{r}$, we have $\norm{G_{i-1}^k(T,x)} \le (K+Kc_{i-1}(k))\norm{x}\le 2Kr$. Hence
\begin{align*}
&\sup_{s \in [iT,iT+T]}\norm{P^k z(s-iT,G_{i-1}^k(T,x))-\phi(s-iT,G_{i-1}^k(T,x))}\\
&\le \sup_{s \in [iT,iT+T]}d^k\norm{\phi(s-iT,G_{i-1}^k(T,x))} \le K d^k\norm{G_{i-1}^k(T,x)}\\
&\le K d^k (K+c_{i-1}(k))\norm{x}.
\end{align*}

Now, if $x\in B_{r}$, $\norm{\phi(s,x)}\le Kr$ and
since $\norm{G_{i-1}^k(T,x)} \le 2Kr$ and $f$ is Lipschitz on $B_{4Kr}$, it is well-known that
\begin{align*}
&\sup_{s \in [iT,iT+T]} \norm{\phi(s-iT,G_{i-1}^k(T,x))-\phi(s-iT,\phi(iT,x))}\\
&\le \sup_{s \in [iT,iT+T]}e^{L(s-iT)}\norm{G_{i-1}^k(T,x)-\phi(iT,x)}\le e^{TL} c_{i-1}(k) \norm{x}
\end{align*}

Combining, we conclude that
\begin{align*}
&\sup_{s \in [iT,iT+T]}\norm{G^k_i(s-iT,x)-\phi(s,x)}\\
&\le e^{TL}  c_{i-1}(k) \norm{x}+K d^k (K+c_{i-1}(k))\norm{x}\\
&= ( (e^{d} + K d^k) c_{i-1}(k) +K^2 d^k )\norm{x} = c_i(k) \norm{x}.
\end{align*}
Since $c_{i}(k) \le c(k)$, and $\delta < NT$, by induction, we conclude that
\begin{equation*}
\sup_{s\in [0,\delta]}\norm{G^k(s,x)-\phi(s,x)}\le  c(k)\norm{x}.
\end{equation*}

\end{proof}

\subsection{Derivative Inequality Lemma} \label{se:derivativeinequalitylemma}

In this critical lemma, we show that the Picard iteration approximately retains the differentiability properties of the solution map. The proof is based on induction, with a key step based on an approach in~\cite{khalil_2002} (Proof of Thm 4.14). This lemma is then adapted to the extended Picard iteration introduced in the previous section.

\begin{lem}\label{lem:der_ineq_Pk}
Suppose that the conditions of Lemma~\ref{lem:contraction_Pk} are satisfied. Then for any $x \in B_r$ and any $k\ge0$,
\begin{align*}
& \sup_{t \in [0,T]}\left\|\frac{\partial}{\partial x}(P^k z)(t,x)^T f(x) -\frac{\partial}{\partial t}(P^k z)(t,x) \right\| \le \frac{(TL)^k}{ T} \norm{x}
\end{align*}
\end{lem}
\begin{proof}
Begin with the identity for $k\ge1$
\begin{align*}
(P^k z)(t,x)&=x + \int_{0}^t f((P^{k-1} z)(s,x))ds\\
&=x + \int_{-t}^0 f((P^{k-1} z)(s+t,x))ds.
\end{align*}
Then, by differentiating the right-hand side, we get
\begin{align*}
&\frac{\partial}{\partial t}(P^k z)(t,x)\\
& =f((P^{k-1} z)(0,x))  \\
&\qquad\qquad +  \int_{-t}^0 \nabla f((P^{k-1} z)(s+t,x))^T\frac{\partial}{\partial 1}(P^{k-1} z)(s+t,x)ds\\
& =f((P^{k-1} z)(0,x)) + \int_{0}^t \nabla f((P^{k-1} z)(s,x))^T\frac{\partial}{\partial s}(P^{k-1} z)(s,x)ds\\
& =f(x)+\int_{0}^t \nabla f((P^{k-1} z)(s,x))^T\frac{\partial}{\partial s}(P^{k-1} z)(s,x)ds,
\end{align*}
where $\frac{\partial}{\partial i}f$ denotes partial differentiation of $f$ with respect to its $i$th variable and
\begin{align*}
&\frac{\partial}{\partial x}(P^k z)(t,x) =I + \int_{0}^t \nabla f((P^{k-1} z)(s,x))^T \frac{\partial}{\partial x}(P^{k-1} z)(s,x)ds.
\end{align*}
Now define for $k\ge1$,
\begin{align*}
y_k(t,x)&:=\frac{\partial}{\partial x}(P^k z)(t,x)^T f(x)-\frac{\partial}{\partial t}(P^k z)(t,x).
\end{align*}
For $k\ge2$, we have
\begin{align*}
&y_k(t,x):=\frac{\partial}{\partial x}(P^k z)(t,x)^T f(x)-\frac{\partial}{\partial t}(P^k z)(t,x)\\
&=-\int_{0}^t \nabla f((P^{k-1} z)(s,x))^T\frac{\partial}{\partial t}(P^{k-1} z)(s,x)ds \\ & \quad + \int_{0}^t \nabla f((P^{k-1} z)(s,x))^T \frac{\partial}{\partial x}(P^{k-1} z)(s,x) f(x)ds\\
&=\int_{0}^t \nabla f((P^{k-1} z)(s,x))^T \cdot \\
& \qquad \qquad \qquad \left[ \frac{\partial}{\partial x}(P^{k-1} z)(s,x) f(x)-\frac{\partial}{\partial s}(P^{k-1} z)(s,x) \right]ds\\
&=\int_{0}^t \nabla f((P^{k-1} z)(s,x))^T y_{k-1}(s,x) ds.
\end{align*}
This means that since $(P^{k-1} z)(t,x) \in B_{2r}$, by induction
\begin{align*}
\sup_{[0,T]}\norm{y_k(t)} &\le T \sup_{t\in[0,T]}\norm{\nabla f((P^{k-1} z)(t,x))} \sup_{t\in[0,T]} \norm{y_{k-1}(t,x)}\\
&\le T L \sup_{t\in[0,T]} \norm{y_{k-1}(t,x)} \le (T L)^{(k-1)} \sup_{t\in[0,T]} \norm{y_{1}(t,x)}
\end{align*}
For $k=1$, $(Pz)(t,x) = x$, so $y_1(t)=f(x)$ and $\sup_{[0,T]}\norm{y_1(t)} \le L\norm{x}$. Thus
\begin{equation*}
\sup_{t\in[0,T]}\norm{y_k(t)} \le \frac{(TL)^k}{T} \norm{x}.
\end{equation*}
\end{proof}

We now adapt this lemma to the extended Picard iteration.

\begin{lem}\label{lem:der_ineq_Gk}
Suppose that the conditions of Lemma~\ref{lem:converge_Gk} are satisfied. Then for any $x \in B_{r}$,
\begin{align*}
&\sup_{t \in [0,T]}\left\|\frac{\partial}{\partial x}G^k (t,x)^T f(x) - \frac{\partial}{\partial t}G^k (t,x)\right\| \le \frac{(TL)^k}{T} (K+ c(k)) \norm{x}
\end{align*}
\end{lem}
\begin{proof}
Recall that
\begin{equation*}
G^k(t,x):=G_i(t-iT,x) \quad \forall \quad t\in [iT,iT+T] \quad \text{and } i=1,\cdots,\infty. 
\end{equation*}
and $G^k_{i+1}(t,x)=P^k z(t,G^k_i(T,x))$ where $z=0$. Then for $t \in [iT,iT+T]$,
\begin{align*}
&\left\|\frac{\partial}{\partial x}G^k (t,x)^T f(x)-\frac{\partial}{\partial t}G^k (t,x)\right\|\\
 &= \left\|\frac{\partial}{\partial x}G^k (t-iT,x)^T f(x)-\frac{\partial}{\partial t}G_i^k (t-iT,x)\right\|\\
&  = \left\|-\frac{\partial}{\partial t}P^k z(t-iT,G^k_i(T,x))+\frac{\partial}{\partial x}P^k (t-iT,G^k_i(T,x))^T f(x)\right\|\\
& \qquad  \le \frac{(TL)^k}{ T} \norm{G^k_i(T,x)}
\end{align*}
As was shown in the proof of Lemma~\ref{lem:converge_Gk}, $\norm{G^k_i(T,x)}\le (K+c_i(k))\norm{x}$. Thus for $t \in [iT,iT+T]$,
\begin{align*}
&\left\|\frac{\partial}{\partial x}G^k (t,x)^T f(x)-\frac{\partial}{\partial t}G^k (t,x) \right\| \le \frac{(TL)^k}{ T} (K+c_i(k))\norm{x}
\end{align*}
Since the $c_i$ are non-decreasing,
\begin{align*}
& \sup_{t \in [0,\delta]}\left\|\frac{\partial}{\partial x}G^k (t,x)^T f(x)-\frac{\partial}{\partial t}G^k (t,x)\right\| \le \frac{(TL)^k}{T} (K+ c(k)) \norm{x}.
\end{align*}
\end{proof}

\section{Main Result - A Converse SOS Lyapunov Function} \label{se:mainresult}

In this section, we combine the previous results to obtain a converse Lyapunov function which is also a sum-of-squares polynomial. Specifically, we use a standard form of converse Lyapunov function and substitute our extended Picard iteration for the solution map. Consider the system
\begin{align}
\dot x (t)=f(x(t)),\qquad x(0)=x_0.\label{eqn:sys}
\end{align}

\begin{theorem}\label{thm:main}
Suppose that $f$ is polynomial of degree $q$ and that system~\eqref{eqn:sys} is exponentially stable on $M$ with
\begin{equation*}
\norm{x(t)}\le K \norm{x(0)}e^{-\lambda t},
\end{equation*}
where $M$ is a bounded nonempty region of radius $r$. Then there exist $\alpha,\beta,\gamma>0$ and a sum-of-squares polynomial $V(x)$ such that for any $x \in M$,
\begin{align}
\alpha \norm{x}^{2} \le V(x)\le \beta \norm{x}^{2} \label{eqn:V_cond1} \\
\nabla V(x)^T f(x)\le  - \gamma \norm{x}^{2}.\label{eqn:V_cond2}
\end{align}
Further, the degree of $V$ will be less than $2q^{ (N k-1)}$, where $k(L,\lambda,K)$ is any integer such that
$c(k)<K$,
\begin{align}
&c(k)^2 + \frac{\log 2K^2}{2 \lambda} K \frac{(TL)^k}{T} (1+c(k)) (K+c(k)) < \frac{1}{2}, \label{eqn:k_cond1} \\
&c(k)^{2} <  \frac{ \lambda}{KL\log 2K^2}(1-(2K^2)^{-\frac{L}{\lambda}}).\label{eqn:k_cond2}
\end{align}
where $c(k)$ is defined as
\begin{equation}
c(k)=\sum_{i=0}^{N-1} \left(e^{TL} + K (TL)^k\right)^i K^2 (TL)^k,\label{eqn:k_cond3}
\end{equation}
  and $N(L,\lambda,K)$ is any integer such that $NT>\frac{\log 2K^2}{2 \lambda}$ and $T<\frac{1}{2L}$ for some $T$ and where $L$ is a Lipschitz bound on $f$ on $B_{4Kr}$.
\end{theorem}


\begin{proof}
Define $\delta = \frac{\log 2K^2}{2 \lambda}$ and $d=TL$. By assumption $N>\frac{\delta}{T}$. Next, we note that since stability implies $f(0)=0$, $f$ is bounded on any $B_r$ with bound $Q=Lr$. Thus for $B_{4Kr}$, we have the bound $Q=4KrL$. By assumption, $T < \frac{1}{2L} = \frac{2Kr}{4KrL} = \frac{2Kr}{Q}$. Therefore, if $k$ is defined as above, the conditions of Lemma~\ref{lem:converge_Gk} are satisfied.  Define $G^k$ as in Lemma~\ref{lem:converge_Gk}. By Lemma~\ref{lem:converge_Gk}, if $k$ is defined as above, $\norm{G^k(s,x)-\phi(s,x)}\le c(k)\norm{\phi(s,x)}$ on $s\in[0,\delta]$ and $x \in B_{r}$.


We propose the following Lyapunov functions, indexed by $k$.
\begin{equation*}
V_k(x):=\int_{0}^{\delta}G^k(s,x)^T G^k(s,x) ds
\end{equation*}
We will show that for any $k$ which satisfies Inequalities~\eqref{eqn:k_cond1},~\eqref{eqn:k_cond2} and~\eqref{eqn:k_cond3}, then if we define $V(x) = V_k(x)$, we have that $V$ satisfies the Lyapunov Inequalities~\eqref{eqn:V_cond1} and~\eqref{eqn:V_cond2} and has degree less than $2q^{(N k-1)}$. The proof is divided into four parts:

\noindent \textbf{Upper and Lower Bounded:}
To prove that $V_k$ is a valid Lyapunov function, first consider upper boundedness. If $x\in B$ and $s \in [0,\delta]$. Then
\begin{align*}
\norm{G^k(s,x)}^2&=\norm{\phi(s,x)+\left[G^k(s,x)^T-\phi(s,x)\right]}^2\\
&\le\norm{\phi(s,x)}^2+\norm{\left[G^k(s,x)^T-\phi(s,x)\right]}^2
\end{align*}
As per Lemma~\ref{lem:converge_Gk}, $\norm{G^k(s,x)-\phi(s,x)}\le c(k) \norm{\phi(s,x)} \le K c(k) \norm{x}$. From stability we have $\norm{\phi(s,x)} \le K \norm{x}$.  Hence,
\begin{equation*}
V_k(x)=\int_{0}^{\delta} \norm{G^k (s,x)}^2 ds \le \delta K^2 \left(1 + c(k)^{2} \right)\norm{x}^2.
\end{equation*}
Therefore the upper boundedness condition is satisfied for any $k\ge 0$ with $\beta = \delta K^2 (1+c(k)^2)>0$.

Next we consider the strict positivity condition. First we note
\begin{align*}
\norm{\phi(s,x)}^2&=\norm{G^k(s,x)+\left[\phi(s,x)-G^k(s,x)\right]}^2\\
&\le\norm{G^k (s,x)}^2+\norm{\phi(s,x)-G^k (s,x)}^2
\end{align*}
which implies
\begin{equation*}
\norm{G^k (s,x)}^2 \ge \norm{\phi(s,x)}^2 -\norm{\phi(s,x)-G^k (s,x)}^2
\end{equation*}
By Lipschitz continuity of $f$, $\norm{\phi(s,x)}^2 \ge e^{-2Ls} \norm{x}^2$ and \\$
\norm{G^k (s,x)-\phi(s,x)}\le K c(k)\norm{x}$.
Thus
\begin{equation*}
V_k(x)=\int_{0}^{\delta} \norm{G^k (s,x)}^2 ds \ge \left( \frac{1}{2L}(1-e^{-2L\delta}) - \delta K c(k)^{2} \right)\norm{x}^2.
\end{equation*}
Therefore for $k$ as defined previously, $\frac{1}{2L}(1-e^{-2L\delta}) - \delta K c(k)^{2} >0$ and so the positivity condition holds for some $\alpha>0$.

\noindent \textbf{Negativity of the Derivative:}
Next, we prove the derivative condition. Recall
\begin{align*}
 V_k(x)&:=\int_{0}^{\delta}G^k(s,x)^T G^k(s,x) ds\\
& = \int_{t}^{t+\delta}G^k(s-t,x)^T G^k(s-t,x) ds
\end{align*}
then since $\nabla V(x(t))^T f(x(t))= \frac{d}{dt}V(x(t))$, we have by the Leibnitz rule for differentiation of integrals,
\begin{align*}
&\frac{d }{d t}V_k(x(t))=\left[G^k(\delta,x(t))^T G^k (\delta,x(t))\right] -\left[G^k(0,x(t))^T G^k(0,x(t))\right]  \\%
&\quad - \int_{t}^{t+\delta}  2 G^k(s-t,x(t))^T \frac{\partial}{\partial 1}G^k(s-t,x(t)) ds \\
&\quad + \int_{t}^{t+\delta} 2 G^k(s-t,x(t))^T \frac{\partial }{\partial 2}G^k(s-t,x(t))f(x(t)) ds\\%
&=\norm{G^k (\delta,x(t))}^2 - \norm{x(t)}^2  \\
&\quad \hspace{-1mm}+ \hspace{-.5mm}\hspace{-.5mm}\int_{0}^{\delta}\hspace{-.5mm}2G^k(s,x(t))^T \left[  \frac{\partial }{\partial 2}G^k(s,x(t))f(x(t))-\frac{\partial}{\partial s}G^k (s,x(t))\right] ds
\end{align*}
where recall $\frac{\partial}{\partial i}f$ denotes partial differentiation of $f$ with respect to its $i$th variable. As per Lemma~\ref{lem:der_ineq_Gk}, we have
\begin{align*}
&\norm{ \frac{\partial }{\partial 2}G^k(t,x(t))^T \hspace{-1mm}f(x(t))-\hspace{-.5mm}\frac{\partial}{\partial 1}G^k (t,x(t))\hspace{-.5mm}} \le \frac{d^k}{T} (K\hspace{-.5mm}+\hspace{-.5mm}c(k)) \norm{x(t)}
\end{align*}
and as previously noted $\norm{ G^k (\delta,x(t))}^2 \le (K^2e^{-2\lambda (s-t)}+c(k)^2)\norm{x(t)}^2$. Also, $\norm{ G^k (s,x(t))}\le K(1+c(k))\norm{x(t)}$. We conclude that
\begin{align*}
&\frac{d }{d t}V_k(x(t))\le(K^2 e^{-2\lambda \delta}+c(k)^2)\norm{x(t)}^2 - \norm{x(t)}^2 \\
&\qquad \qquad \qquad \qquad + 2 \delta \frac{d^k}{T} K (1 + c(k)) (K+c(k)) \norm{x(t)}^2 \\
&\le\hspace{-.5mm} \hspace{-.5mm}\left(K^2 e^{-2\lambda \delta}\hspace{-.5mm}+c(k)^2\hspace{-.5mm}\hspace{-.5mm}-\hspace{-.5mm}1 + 2 \delta K\frac{d^k}{T} (1 + c(k)) (K+c(k)) \hspace{-.5mm}\hspace{-.5mm}\right)\hspace{-.5mm}\hspace{-.5mm}\norm{x(t)}^2\hspace{-.5mm}\hspace{-.5mm}.
\end{align*}
Therefore, we have strict negativity of the derivative since
\begin{align*}
&K^2e^{-2\lambda \delta}+c(k)^2 + 2\delta\frac{d^k}{T} (1+c(k)) (K+c(k)) \\
&\quad = \frac{1}{2}+c(k)^2 + 2\frac{K\log 2K^2}{2 \lambda} \frac{d^k}{T} (1+c(k)) (K+c(k))<1
\end{align*}
Thus $\frac{d }{d t}V_k(x(t))\le - \gamma\norm{x(t)}^2$ for some $\gamma>0$.



\noindent \textbf{Sum of Squares:}
Since $f$ is polynomial and $z$ is trivially polynomial, $(P^k z)(s,x)$ is a polynomial in $x$ and $s$. Therefore, $V_k(x)$ is a polynomial for any $k>0$. To show that $V_k$ is sum-of-squares, we first rewrite the function
\begin{align*}
V_k(x)&=\sum_{i=1}^N \int_{iT-T}^{iT}  \left[G_i^k (s-iT,x)^T G_i^k(s-iT,x)\right] ds.
\end{align*}
Since $G_i^k z$ is a polynomial in all of its arguments, $G_i^k (s-iT,x)^T G_i^k(s-iT,x)$ is sum-of-squares. It can therefore be represented as $R_i(x)^T Z_i(s)^T Z_i(s)R_i(x)$ for some polynomial vector $R_i$ and matrix of monomial bases $Z_i$. Then
\begin{align*}
V_k(x)&= \sum_{i=1}^N R_i(x)^T \int_{iT-T}^{iT} Z_i(s)^T Z_i(s)ds R_i(x) = \sum_{i=1}^N R_i(x)^T M_i R_i(x)
\end{align*}
Where $M_i=\int_{iT-T}^{iT} Z_i(s)^T Z_i(s)ds \ge 0$ is a constant matrix. This proves that $V_k$ is sum-of-squares since it is a sum of sums-of-squares.

We conclude that $V=V_k$ satisfies the conditions of the theorem for any $k$ which satisfies Inequalities~\eqref{eqn:k_cond1} and~\eqref{eqn:k_cond2}.
\noindent \textbf{Degree Bound:} Given a $k$ which satisfies the inequality conditions on $c(k)$, we consider the resulting degree of $G^k$, and hence, of $V_k$. If $f$ is a polynomial of degree $q$, and $y$ is a polynomial of degree $d$ in $x$, then $P y$ will be a polynomial of degree $\max \{1,dq\}$ in $x$. Thus since $z=0$, the degree of $P^k z$ will be $q^{k-1}$. If $N>1$, then the degree of $G_i^k$ will be $q^{Nk-1}$. Thus the maximum degree of the Lyapunov function is $2q^{(Nk-1)}$.

\end{proof}

In the proof of Theorem~\ref{thm:main}, the integration interval, $\delta$ was chosen such that the conditions will always be feasible for some $k>0$. However, this choice may not be optimal. Numerical experimentation has shown us that a better degree bound may be obtained by varying this parameter in the proof. However, the given value is one which we have found to work well in the vast majority of cases.

We conclude this section by commenting on the form of the converse Lyapunov function,
\begin{equation*}
V_k(x):=\int_{0}^{\delta}G^k(s,x)^T G^k(s,x) ds.
\end{equation*}
Our Lyapunov function is defined using an approximation of the solution map. A dual approach to solution of the Hamilton-Jacobi-Bellmand Equation was taken in~\cite{lasserre_2008} using occupation measures instead of Picard iteration. Indeed, the dual space of the Sum of Squares Lyapunov functions
can be understood in terms of moments of such occupation measures~\cite{peyrl_2007}.

As a final note, the proof of Theorem~\ref{thm:main} also holds for time-varying systems. Indeed the original proof was for this case. However, because Sum-of-Squares is rarely used for time-varying systems, the result has been simplified to improve clarity of presentation.


\subsection{Numerical Illustration}~\label{se:numericalexample}
To illustrate the degree bound and hence the complexity of analyzing a nonlinear system, we plot the degree bound versus the exponential convergence rate of the system. For given parameters, this bound is obtained by numerically searching for the smallest $k$ which satisfies the conditions of Theorem~\ref{thm:main}. The convergence rate parameter can be viewed as a metric for the accuracy of the sum-of-squares approach: suppose we have a degree bound as a function of convergence rate, $d(\gamma)$. If it is not possible to find a sum-of-squares Lyapunov function of degree $d(\gamma)$ proving stability, then we know that the convergence rate of the system must be less than $\gamma$.

\begin{figure}
 \begin{tabular}{c}
  \includegraphics[scale=.42]{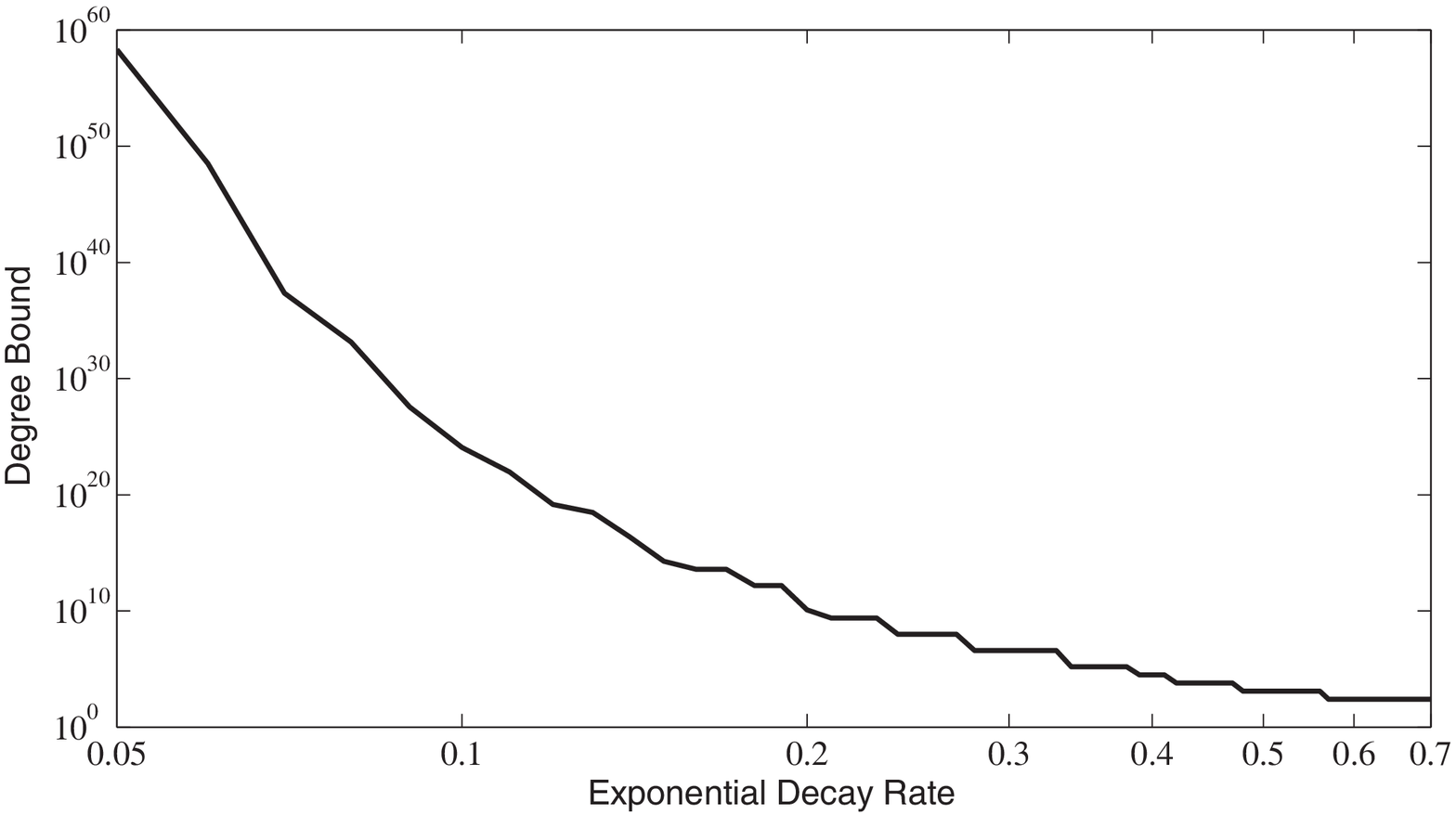} \\
  \includegraphics[scale=.42]{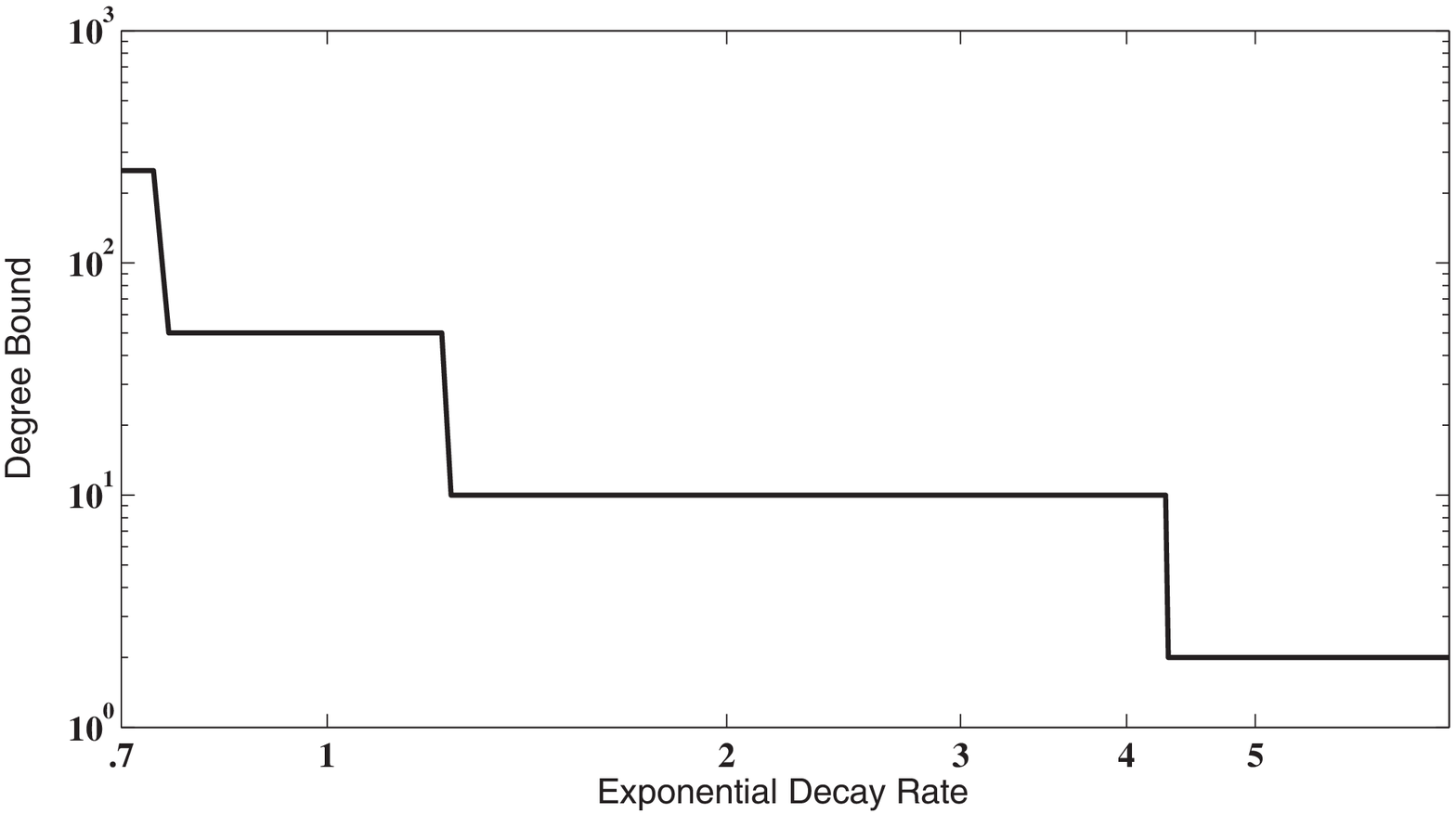}
 \end{tabular}
 \caption{Degree bound vs. Exponential Convergence Rate for $K=1.2$, $r=L=1$, $q=5$. Domains $\lambda < .7$ and $\lambda >.7$ are plotted separately for clarity.}\label{f:test}
\end{figure}

As can be seen, as the convergence rate increases, the degree bound decreases super-exponentially, so that at $\gamma=2.4$, only a quadratic Lyapunov function is required to prove stability. For cases where high accuracy is required, the degree bound increases quickly; scaling approximately as $e^{\frac{1}{\gamma}}$. To reduce the complexity of the problem, in come cases less conservative bounds on the degree can be found by considering the monomial terms in the vector field. If the complexity is still unacceptably high, then one can consider the use of parallel computing: unlike single-core processing, parallel computing power continues to increase exponentially. For a discussion on using parallel computing to solve polynomial optimization problems, we refer to~\cite{PeeP10}.

\section{Quadratic Lyapunov Functions}\label{se:quadratic}
In this section, we briefly explore the implications of our result for the existence of quadratic Lyapunov functions proving exponential stability of nonlinear systems. Specifically, we look at when the theorem predicts the existence of a degree bound of 2. We first note that when the vector field is linear, then $q=1$, which implies that $2q^{(Nk-1)}=2$ independent of $N$ and $k$. Recall $N$ is the number of Picard iterations, $k$ is the number of extensions and $q$ is the degree of the polynomial vector field, $f$. Hence an exponentially stable linear system has a quadratic Lyapunov function - which is not surprising.

Instead we consider the case when $q \neq 1$. In this case, for a quadratic Lyapunov function, we require $N=k=1$ - a single Picard iteration and no extensions. By examining the proof of Theorem~\ref{thm:main}, we see that if the conditions of the theorem are satisfied with $N=k=1$ then $V(x)=x^T x$ is a Lyapunov function which establishes exponential stability of the system. Since this is perhaps the most commonly used form of Lyapunov function, it is worth considering how conservative it is when applied to nonlinear systems of the form
\[
\dot x(t) = f(x(t)).
\]
In the following corollary we give sufficient conditions on the vector field and decay rate for the Lyapunov function $x^T x$ to prove exponential stability.
\begin{cor}
Suppose that system~\eqref{eqn:sys} is exponentially stable with
\begin{equation*}
\norm{x(t)}\le K \norm{x(0)}e^{-\lambda t}
\end{equation*}
for some $\lambda>0$, $K\ge1$ and for any $x(0)\in M$, where $M$ is a bounded nonempty region of radius $r$. Let $L$ be a Lipschitz bound for $f$ on $B_{4Kr}$.
Suppose that there exists some $\frac{1}{2L}> \delta >0$ such that
\[
K^2 e^{-2 \lambda \delta}+c_1^2+2 K \delta L (1+c_1)(K+c_1) < 1
\]
and $K \delta L <1 $, where $c_1 = K^2 \delta L$.
Let $V(x)=x^Tx$. Then for any $x \in M$,
\begin{align*}
\dot V(x) = \nabla x^T f(x) \le  - \beta \norm{x}^{2}.
\end{align*}
for some $\beta>0$.
\end{cor}
\begin{proof}
We reconsider the proof of Theorem~\ref{thm:main}. This time, we set $N=k=1$ and $T=\delta$ and determine if there exists a $\delta=T<\frac{1}{2L}$ which satisfies the upper-boundedness, lower-boundedness and derivative conditions. Because $V(x) = \delta x^T x$, the upper and lower boundedness conditions are immediately satisfied. The derivative negativity condition is
\[
K^2 e^{-2\lambda \delta}+c(1)^2 + 2K \delta L (1+c(1)) (K+c(1)) < 1
\]
where $c(1)=c_1 = K^2 \delta L$. This is satisfied by the statement of the theorem.
\end{proof}

\begin{figure}
\includegraphics[scale=.45]{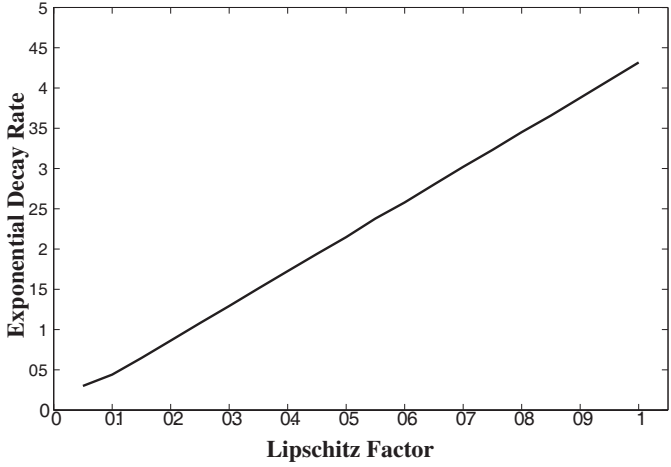}
 \caption{Required decay rate for a quadratic Lyapunov function vs. Lipschitz bound for $K=1.2$}\label{f:quadratic}
\end{figure}

Note that neither the size of the region we consider nor the degree of the vector field plays any role in determining the degree bound.
To illustrate the conditions for existence of a quadratic Lyapunov function, we plot the required decay rate vs. the Lipschitz continuity factor in Figure~\ref{f:quadratic} for $K=1.2$. This plot shows that as the Lipschitz continuity of the vector field increases (and the field becomes less smooth), the conservatism of using the quadratic Lyapunov function $x^T x$ increases.

\section{Implications for Sum-of-Squares Programming} \label{se:benefits}

In this section we consider the implications that the above results have on Sum of Squares programming.

\subsection{Bounding the number of decision variables}

Because the set of continuously differentiable functions is an infinite-dimensional vector space, the general problem of finding a Lyapunov function is an infinite-dimensional feasibility problem. However, the set of sum-of-squares Lyapunov functions with bounded degree is finite-dimensional. The most significant implication of our theorem is a bound on the number of variables in the problem of determining stability of a nonlinear vector field. The nonlinear stability problem can now be expressed as a feasibility problem of the following form.
\begin{theorem}\label{thm:optim}
For a given $\lambda$, let $2d$ be the degree bound associated with Theorem~\ref{thm:main} and define $N= \frac{(n+d)!}{n! d!}$. If System~\eqref{eqn:sys} is exponentially stable on $M$ with decay rate $\lambda$ or greater, the following is feasible for some $\alpha,\beta,\gamma>0$.
\begin{align*}
\text{\textbf{Find:} } P \in \S^N:\\
P&\ge0\\
\alpha \norm{x}^2 \le Z(x)^T P Z(x) &\le \beta \norm{x}^2 \qquad \text{for all } x \in M\\
\nabla \left(Z(x)^T P Z(x)\right)^T f(x) &\le -\gamma \norm{x}^2\qquad \text{for all } x \in M
\end{align*}
where $Z(x)$ be the vector of monomials in $x$ of degree $d$ or less.
\end{theorem}
\begin{proof}
The proof follows immediately from the fact that a polynomial $V$ of degree $2d$ is SOS if and only if there exists a $P\ge0$ such that $V(x) = Z(x)^T P Z(x)$.
\end{proof}
Our condition bounds the number of variables in the feasibility problem associated with Theorem~\ref{thm:optim}. If $M$ is semialgebraic, then the conditions in Theorem~\ref{thm:optim} can be enforced using sum-of-squares and the Positivstellensatz~\cite{putinar_1993}. The complexity of solving the optimization problem will depend on the complexity of the Positivstellensatz test. If positivity on a semialgebraic set is decidable, as indicated in~\cite{nie_2007}, this implies the question of exponential stability on a bounded set is decidable.

\subsection{Local Positivity}

Another implication of our result is that it reduces the complexity of enforcing the positivity constraint. As discussed in Section~\ref{sec:SOS}, semidefinite programming is used to optimize over the cone of sums-of-squares of polynomials. There are several different ways the stability conditions can be enforced. For example, we have the following theorem.
\begin{theorem}\label{thm:straw}
Suppose there exist polynomial $V$ and sum-of-squares polynomials $s_1,s_2,s_3$ and $s_4$ such that the following conditions are satisfied for $\alpha,\gamma>0$.
\begin{align*}
V(x) - \alpha \norm{x}^2&= s_1(x) + g(x)s_2(x)\\
-\nabla V(s)^T f(x) -\gamma \norm{x}^2 &= s_3(x) + g(x)s_4(x)
\end{align*}
Then we have exponential stability of System~\eqref{eqn:sys} on $\left\{ x \;:\; \{y: V(y) \le V(x)\}\subset U\right\}$.
\end{theorem}

The complexity of the conditions associated with Theorem~\ref{thm:straw} is determined by the four sum-of-squares variables, $s_i$. Theorem~\ref{thm:straw} uses the Positivstellensatz multipliers $s_2$ and $s_4$ to ensure that the Lyapunov function need only be positive and decreasing on the region $X=\{x: g(x) \ge 0\}$. However, as we now know that the Lyapunov function can be assumed SOS, we can eliminate the multiplier $s_2$, reducing complexity of the problem.

\begin{theorem}\label{thm:new_stability}
Suppose there exist polynomial $V$ and sum-of-squares polynomials $s_1$, $s_2$ and $s_3$ such that the following conditions are satisfied for $\alpha,\gamma>0$.
\begin{align*}
V(x) - \alpha \norm{x}^2&= s_1(x)\\
-\nabla (V(x)+\alpha \norm{x}^2)^T f(x) - \gamma \norm{x}^2 &= s_2(x) + g(x)s_3(x)
\end{align*}
Then we have exponential stability of System~\eqref{eqn:sys} for any $x(0)$ such that $\{y: V(y) \le V(x(0))\}\subset X$ where $X:=\{x:g(x)\ge 0\}$.
\end{theorem}

This simplification reduces the size of the SOS variables by 25\% (from 4 to 3). If the semialgebraic set $X$ is defined using several polynomials (e.g. a hypercube), then the reduction in the number of variables can approach 50\% . SDP solvers are typically of complexity $O(n^6)$, where $n$ is the dimension of the symmetric matrix variable. In the above example we reduced $n=4N$ to $n=3N$. Thus this simplification can potentially decrease computation by a factor of 82\%.

\section{Numerical Example} \label{se:example}

\begin{figure}
\includegraphics[scale=.62]{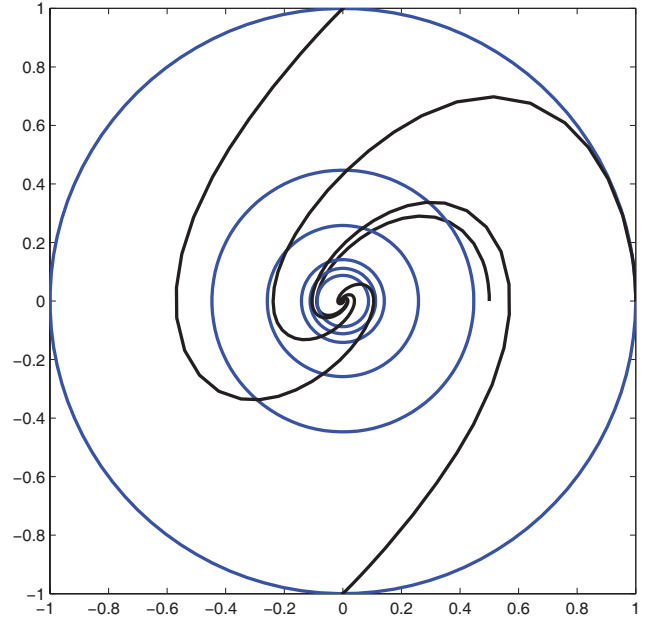}\vspace{-2mm}
 \caption{Plot of trajectories of the Van-der-Pol Oscillator. We estimate the overshoot parameter as $K\cong1$}\label{f:estim_K}
\end{figure}

\begin{figure}
 \includegraphics[scale=.5]{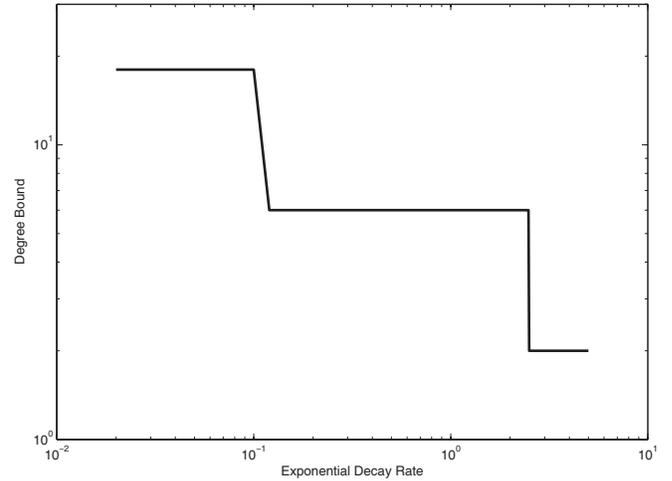}
 \caption{Degree Bound for the Van-der-Pol Oscillator as a Function of Decay Rate}\label{f:dplot_vanderpol}
\end{figure}

\begin{figure}
\includegraphics[scale=.5]{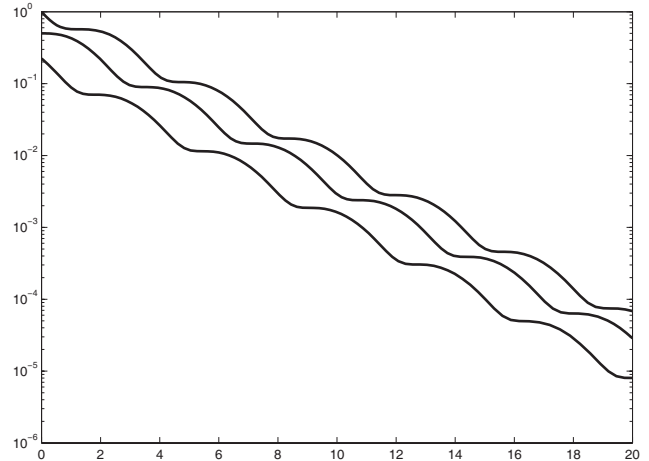}
 \caption{A semi-log plot of $\norm{x}$ for three trajectories. We estimate $\lambda=.542$ for the Van-der-Pol oscillator}\label{f:estim_decay}
\end{figure}

In this section, we use the Van-der-Pol oscillator to illustrate how the degree bound influences the accuracy of the stability test. The zero equilibrium point of the Van-der-Pol oscillator is unstable. In reverse-time, however, this equilibrium is stable with a domain of attraction bounded by the well-known forward-time limit-cycle. The reverse-time dynamics are as follows.
\begin{align*}
\dot x_1(t) &= -x_2(t)\\
\dot x_2(t) &= -\mu (1-x_1(t)^2)x_2(t) +x_1(t)
\end{align*}
For simplicity, we choose $\mu=1$. On a ball of radius $r$, the Lipschitz constant can be found from $L = \sup_{x \in B_r} \norm{D f(x)}$, where $\norm{\cdot}$ is the maximum singular value norm. We find a Lipschitz constant for the Van-der-Pol oscillator on radius $r=1$ to be $2.1$. Numerical simulations indicate $K \cong 1$, as illustrated in Figure~\ref{f:estim_K}.
Given these parameters, the degree bound plot is illustrated in Figure~\ref{f:dplot_vanderpol}. Note that the choice of $K=1$ dramatically improves the degree bound.
Numerical simulation shows the decay rate to be a relatively constant $\lambda = .542$ throughout the unit ball. This is illustrated in Figure~\ref{f:estim_decay}. This gives us an estimate of the degree bound as $d=6$.


%

To find the converse Lyapunov function associated  with this degree bound we construct the Picard iteration.
\begin{align*}
(P z)(t,x) &= x + \int_{0}^t f(0)ds = x.\\
(P^2 z)(t,x) &= x + \int_{0}^t f(P z(s,x))ds\\
&=x + \int_{0}^t f( x )ds=x + f(x) t
\end{align*}
The converse Lyapunov function is
\begin{align*}
V(x) &= \int_{0}^\delta (P^2 z(s,x))^T (P^2 z(s,x)) ds\\
&=\int_{0}^\delta (x + f(x)s)^T (x + f(x)s) ds\\
&=\int_{0}^\delta \bmat{x\\f(x)}^T\bmat{I\\sI}\bmat{I&sI}\bmat{x\\f(x)} ds\\
&= \bmat{x\\f(x)}^T \int_{0}^\delta \bmat{I & sI\\sI & s^2I}ds\bmat{x\\f(x)} \\
&= \bmat{x\\f(x)}^T  \bmat{\delta I & \delta^2/2 I\\ \delta^2/2 I & \delta^3/3 I} \bmat{x\\f(x)}
\end{align*}
If $\delta = T = \frac{1}{2L} = \frac{1}{4}$, for the Van-der-Pol Oscillator, we get the SOS Lyapunov function.
\begin{align*}
&192\cdot V(x) = \bmat{x\\f(x)}^T  \bmat{ 48 I & 6 I\\ 6 I & I}\bmat{x\\f(x)}\\
&= \bmat{x\\f(x)}^T  \bmat{ 6.93 I & 2.45 I\\ 2.45 I & I}^2 \bmat{x\\f(x)}\\
&= \bmat{6.93 x + 2.45 f(x)\\2.45 x + f(x)}^T \bmat{6.93 x+2.45 f(x)\\2.45 x + f(x)}\\
&= \left(6.93 x_1 - 2.45 x_2\right)^2 +\left(2.45 (x_1 +x_1^2 x_2) +4.48 x_2 \right)^2\\
 & \qquad + \left( 2.45 x - x_2\right)^2+ \left( x_1+x_1^2 x_2 + 1.45 x_2\right)^2
\end{align*}

\begin{figure}
\includegraphics[scale=.66]{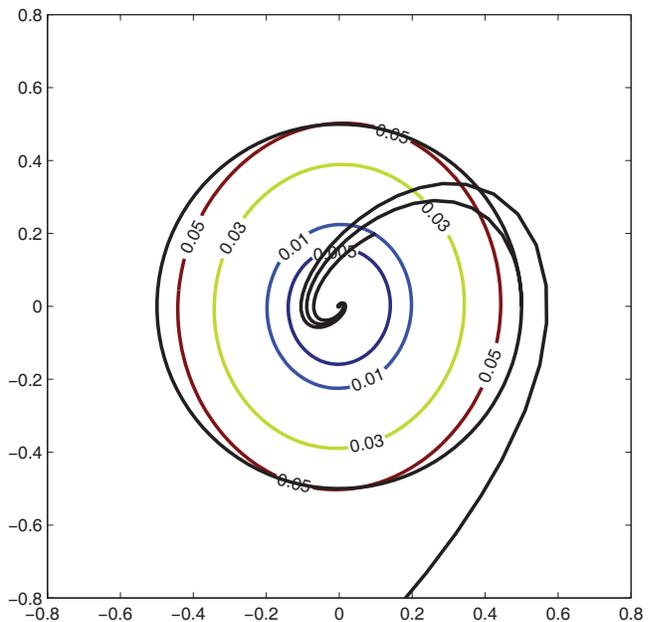}
 \caption{Level Sets of the converse Lyapunov function, with Ball of radius $r=.25$}\label{f:lyap_quadratic}
\end{figure}

As per the previous discussion, we use SOSTOOLS to verify that this Lyapunov function proves stability. Note that we must show the function is decreasing on the ball of radius $r=.25$, as the Lipschitz bound used in the theorem is for the ball of radius $B_{4r}$. We are able to verify that the Lyapunov function is decreasing on the ball of radius $r=.25$. Some level sets of this Lyapunov function are illustrated in Figure~\ref{f:lyap_quadratic}. Through experimentation, we find that when we increase the ball to radius $r=1$, the Lyapunov function is no longer decreasing. We also found that the quadratic Lyapunov function $V(x)=x^T x$ is not decreasing on the ball of radius $r=.25$. Although we believe that our degree bound is somewhat conservative, these results indicate the conservatism is not excessive.

%
%
%

\begin{figure}
\includegraphics[scale=.66]{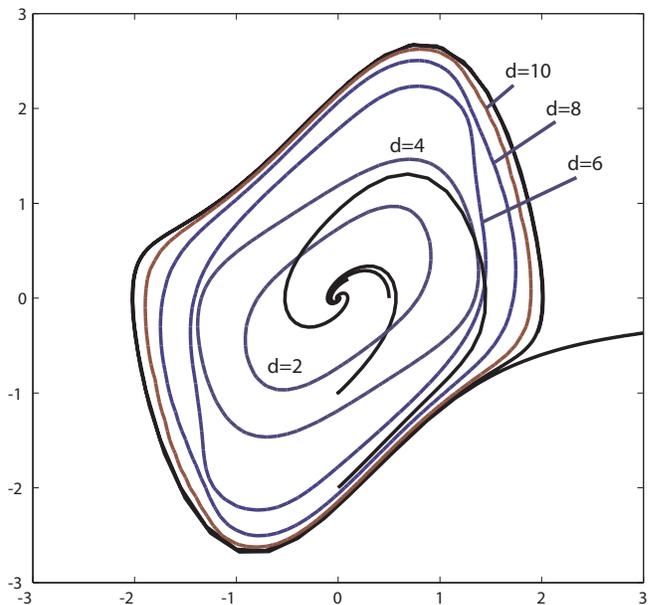}
 \caption{Best Invariant Region vs. Degree Bound with Limit Cycle}\label{f:lyap_composite}
\end{figure}


To explore the limits of the SOS approach, for degree bound 2, 4, 6, 8 and 10, we find the maximum unit ball on which we are able to find a sum-of-squares Lyapunov function. We then use the largest sublevel set of this Lyapunov function on which the trajectories decrease as an estimate for the domain of attraction of the system. These level sets are illustrated in Figure~\ref{f:lyap_composite}. We see that as the degree bound increases, our estimate of the domain of attraction improves.

\section{Conclusion} \label{se:conclusion}

In this paper, we have used the Picard iteration to construct an approximation to the solution map on arbitrarily long intervals. We have used this approximation to prove that exponential stability of a polynomial vector field on a bounded set implies the existence of a Lyapunov function which is a sum-of-squares of polynomials with a bound on the degree. This implies that the question of exponential stability on a bounded set may be decidable. Furthermore, the converse Lyapunov function we have used in this paper is relatively easy to construct given the vector field and may find applications in other areas of control. The main result also holds for time-varying systems.

Recently, there has been interest in using semidefinite programming for the analysis on nonlinear systems using sum-of-squares. This paper clarifies several questions on the application of this method. We now know that exponential stability on a bounded set implies the existence of an SOS Lyapunov function and we know how complex this function may be. It has been recently shown that \textit{globally} asymptotically stable vector fields do not always admit sum-of-squares Lyapunov functions~\cite{ahmadi_2011}. Still unresolved is the question of the existence of polynomial Lyapunov functions for stability of globally \textit{exponentially} stable vector fields.

\bibliographystyle{ieeetr}
\bibliography{weierstrass,delay,lyapunov,biblio}

\begin{thebibliography}{10}

\bibitem{BoyEFB94}
S.~Boyd, L.~{El Ghaoui}, E.~Feron, and V.~Balakrishnan, {\em Linear Matrix
  Inequalities in System and Control Theory}.
\newblock SIAM Studies in Applied Mathematics, 1994.

\bibitem{lasserre_2001}
J.~B. Lasserre, ``Global optimization with polynomials and the problem of
  moments,'' {\em SIAM J. Optim.}, vol.~11, no.~3, pp.~796--817, 2001.

\bibitem{nesterov_2000}
Y.~Nesterov, {\em High Performance Optimization}, vol.~33 of {\em Applied
  Optimization}, ch.~Squared Functional Systems and Optimization Problems.
\newblock Springer, 2000.

\bibitem{parrilo_thesis}
P.~A. Parrilo, {\em Structured Semidefinite Programs and Semialgebraic Geometry
  Methods in Robustness and Optimization}.
\newblock PhD thesis, California Institute of Technology, 2000.

\bibitem{henrion_2005}
D.~Henrion and A.~Garulli, eds., {\em Positive Polynomials in Control},
  vol.~312 of {\em Lecture Notes in Control and Information Science}.
\newblock Springer, 2005.

\bibitem{chesi_2007}
G.~Chesi, ``On the gap between positive polynomials and {SOS} of polynomials,''
  {\em IEEE Transactions on Automatic Control}, vol.~52, pp.~1066--1072, June
  2007.

\bibitem{Che10}
G.~Chesi, ``{LMI} techniques for optimization over polynomials in control: A
  survey,'' {\em IEEE Transactions on Automatic Control}, vol.~55, no.~11,
  pp.~2500--2510, 2010.

\bibitem{prajna_2004}
S.~Prajna, A.~Papachristodoulou, P.~Seiler, and P.~A. Parrilo, ``New
  developments in sum of squares optimization and {SOSTOOLS},'' in {\em
  Proceedings of the American Control Conference}, pp.~5606 -- 5611, 2004.

\bibitem{henrion_2001}
D.~Henrion and J.-B. Lassere, ``{GloptiPoly}: Global optimization over
  polynomials with {MATLAB} and {SeDuMi},'' in {\em IEEE Conference on Decision
  and Control}, pp.~747--752, 2001.

\bibitem{Papachristodoulou_2002}
A.~Papachristodoulou and S.~Prajna, ``On the construction of {L}yapunov
  functions using the sum of squares decomposition,'' in {\em Proceedings IEEE
  Conference on Decision and Control}, pp.~3482 -- 3487, 2002.

\bibitem{wang_thesis}
T.-C. Wang, {\em Polynomial Level-Set Methods for Nonlinear Dynamics and
  Control}.
\newblock PhD thesis, Stanford University, 2007.

\bibitem{tan_2006}
W.~Tan, {\em Nonlinear Control Analysis and Synthesis using Sum-of-Squares
  Programming}.
\newblock PhD thesis, University of California, Berkeley, 2006.

\bibitem{peet_2009a}
M.~M. Peet, ``Exponentially stable nonlinear systems have polynomial {Lyapunov}
  functions on bounded regions,'' {\em IEEE Transactions on Automatic Control},
  vol.~52, pp.~979--987, May 2009.

\bibitem{barbasin_1951}
E.~A. Barbasin, ``The method of sections in the theory of dynamical systems,''
  {\em Rec. Math. (Mat. Sbornik) N. S.}, vol.~29, pp.~233--280, 1951.

\bibitem{malkin_1954}
I.~Malkin, ``On the question of the reciprocal of {Lyapunov's} theorem on
  asymptotic stability,'' {\em Prikl. Mat. Meh.}, vol.~18, pp.~129--138, 1954.

\bibitem{kurzweil_1963}
J.~Kurzweil, ``On the inversion of {Lyapunov's} second theorem on stability of
  motion,'' {\em Amer. Math. Soc. Transl.}, vol.~2, no.~24, pp.~19--77, 1963.
\newblock English Translation. Originally appeared 1956.

\bibitem{massera_1956}
J.~L. Massera, ``Contributions to stability theory,'' {\em Annals of
  Mathematics}, vol.~64, pp.~182--206, July 1956.

\bibitem{wilson_1969}
F.~W. {Wilson Jr.}, ``Smoothing derivatives of functions and applications,''
  {\em Transactions of the American Mathematical Society}, vol.~139,
  pp.~413--428, May 1969.

\bibitem{lin_1996}
Y.~Lin, E.~Sontag, and Y.~Wang, ``A smooth converse {Lyapunov} theorem for
  robust stability,'' {\em Siam J. Control Optim.}, vol.~34, no.~1,
  pp.~124--160, 1996.

\bibitem{lakshmikantam_1992}
V.~Lakshmikantam and A.~A. Martynyuk, ``{Lyapunov's} direct method in stability
  theory (review),'' {\em International Applied Mechanics}, vol.~28,
  pp.~135--144, March 1992.

\bibitem{teel_1999}
A.~R. Teel and L.~Praly, ``Results on converse {Lyapunov} functions from
  class-{KL} estimates,'' pp.~2545--2550, 1999.

\bibitem{hahn_1967}
W.~Hahn, {\em Stability of Motion}.
\newblock Springer-Verlag, 1967.

\bibitem{krasovskii_1963}
N.~N. Krasovskii, {\em Stability of Motion}.
\newblock Stanford University Press, 1963.

\bibitem{arnold_book}
V.~I. Arnol'd, {\em Ordinary Differential Equations}.
\newblock Springer, 2~ed., 2006.
\newblock Translated by {Roger Cook}.

\bibitem{MurK87}
K.~G. Murty and S.~N. Kabadi, ``Some {NP}-complete problems in quadratic and
  nonlinear programming,'' {\em Mathematical Programming}, vol.~39,
  pp.~117--129, 1987.

\bibitem{Sho87}
N.~Z. Shor, ``Class of global minimum bounds of polynomial functions,'' {\em
  Cybernetics}, vol.~23, no.~6, pp.~731--734, 1987.

\bibitem{PowW98}
V.~Powers and T.~W\"ormann, ``An algorithm for sums of squares of real
  polynomials,'' {\em Journal of Pure and Applied Linear Algebra}, vol.~127,
  pp.~99--104, 1998.

\bibitem{Par00}
P.~A. Parrilo, {\em Structured Semidefinite Programs and Semialgebraic Geometry
  Methods in Robustness and Optimization}.
\newblock PhD thesis, Caltech, Pasadena, CA, 2000.
\newblock Available at \verb"http://www.mit.edu/~parrilo/pubs/index.html".

\bibitem{marsden_analysis}
J.~E. Marsden and M.~J. Hoffman, {\em Elementary Classical Analysis}.
\newblock W. H. Feeman and Company, 2nd~ed., 1993.

\bibitem{coddingtion_1955}
E.~A. Coddington and N.~Levinson, {\em Theory of Ordinary Differential
  Equations}.
\newblock McGraww-Hill, 1955.

\bibitem{lindelhof_1894}
E.~Lindel\"{o}f and M.~Picard, ``Sur l'application de la m\'{e}thode des
  approximations successives aux \'{e}quations diff\'{e}rentielles ordinaires
  du premier ordre,'' {\em Comptes rendus hebdomadaires des s\'{e}ances de
  l'Acad\'{e}mie des sciences}, vol.~114, pp.~454--457, 1894.

\bibitem{khalil_2002}
H.~Khalil, {\em Nonlinear Systems}.
\newblock Prentice Hall, third~ed., 2002.

\bibitem{lasserre_2008}
J.~B. Lasserre, D.~Henrion, C.~Prieur, and E.~Trelat, ``Nonlinear optimal
  control via occupation measures and {LMI}-relaxations,'' vol.~47, no.~4,
  pp.~1643--1666, 2008.

\bibitem{peyrl_2007}
H.~Peyrl and P.~A. Parrilo, ``A theorem of the alternative for {SOS Lyapunov}
  functions,'' in {\em Proceedings IEEE Conference on Decision and Control},
  pp.~1687--1692, 2007.

\bibitem{PeeP10}
M.~M. Peet and Y.~V. Peet, ``A parallel-computing solution for optimization of
  polynomials,'' in {\em Proceedings of the American Control Conference},
  pp.~4851 -- 4856, 2010.

\bibitem{putinar_1993}
M.~Putinar, ``Positive polynomials on compact semi-algebraic sets,'' {\em
  Indiana Univ. Math. J.}, vol.~42, no.~3, pp.~969--984, 1993.

\bibitem{nie_2007}
J.~Nie and M.~Schweighofer, ``On the complexity of {Putinar's}
  positivstellensatz,'' {\em Journal of Complexity}, vol.~23, pp.~135--150,
  2007.

\bibitem{ahmadi_2011}
A.~A. Ahmadi, M.~Krstic, and P.~A. Parrilo, ``A globally asymptotically stable
  polynomial vector field with no polynomial {Lyapunov} function,'' in {\em
  Proceedings of the IEEE Conference on Decision and Control}, pp.~7579--7580,
  2011.

\end{thebibliography}
\vspace{-3mm}
\begin{IEEEbiography}{Matthew M. Peet}
 received B.S. degrees in Physics and in Aerospace Engineering from
the University of Texas at Austin in 1999 and the M.S. and Ph.D. in
Aeronautics and Astronautics from Stanford University in 2001 and
2006, respectively. He was a Postdoctoral Fellow at the National
Institute for Research in Computer Science and Control (INRIA) near
Paris, France, from 2006-2008 where he worked in the SISYPHE and BANG groups. He is currently an
Assistant Professor in the Mechanical, Materials, and Aerospace
Engineering Department of the Illinois Institute of Technology and
director of the Cybernetic Systems and Controls Laboratory.
His current research interests are in the role of computation as it
is applied to the understanding and control of complex and
large-scale systems. Applications include fusion energy and immunology.
\end{IEEEbiography}
\vspace{-3mm}
\begin{IEEEbiography}{Antonis Papachristodoulou}
 received an MA/MEng degree in Electrical and
Information Sciences from the University of Cambridge in 2000, as a
member of Robinson College. In 2005 he received a Ph.D. in Control and
Dynamical Systems, with a minor in Aeronautics from the California
Institute of Technology. In 2005 he held a David Crighton Fellowship
at the University of Cambridge and a postdoctoral research associate
position at the California Institute of Technology before joining the
Department of Engineering Science at the University of Oxford, Oxford,
UK in January 2006, where he is now a University Lecturer in Control
Engineering and tutorial fellow at Worcester College. His research
interests include scalable analysis of nonlinear systems using convex
optimization based on Sum of Squares programming, analysis and design
of large-scale networked control systems with communication
constraints and Systems and Synthetic Biology.
\end{IEEEbiography}

\end{document}